\documentclass[11pt]{amsart}
\usepackage[
  textheight=615pt,
  textwidth=360pt,
  centering
]{geometry}
\usepackage{amssymb,bbm,color,xcolor}
\usepackage{amscd}
\usepackage{enumerate}
\usepackage{siunitx}
\usepackage{bm}
\usepackage{ulem}
\usepackage{float, graphicx}
\usepackage{subfigure}
\numberwithin{equation}{section}

\usepackage[pagebackref]{hyperref}

\usepackage{mathtools}
\usepackage[tableposition=top]{caption}
\usepackage{booktabs,dcolumn}

\DeclareFontFamily{OT1}{rsfs}{}
\DeclareFontShape{OT1}{rsfs}{n}{it}{<-> rsfs10}{}
\DeclareMathAlphabet{\mathscr}{OT1}{rsfs}{n}{it}

\theoremstyle{plain}

\newtheorem{thm}{Theorem}[section]
\newtheorem{prop}[thm]{Proposition}
\newtheorem{lem}[thm]{Lemma}

\newtheorem{conjecture}[thm]{Conjecture}

\theoremstyle{definition}

\newtheorem{remark}[thm]{Remark}

\newcommand{\bal}{\[\begin{aligned}}
\newcommand{\eal}{\end{aligned}\]}
\newcommand{\beeq}{\begin{equation}}\newcommand{\eneq}{\end{equation}}
\newcommand{\beq}{\begin{equation*}}\newcommand{\eeq}{\end{equation*}}
\newcommand{\beqa}{\begin{eqnarray*}}\newcommand{\eeqa}{\end{eqnarray*}}
\newcommand{\beeqa}{\begin{eqnarray}}\newcommand{\eneqa}{\end{eqnarray}}
\newcommand{\beitm}{\begin{itemize}}\newcommand{\eitm}{\end{itemize}}

\def\<{\langle}             \def\>{\rangle}
\newcommand{\al}{\alpha}    \newcommand{\be}{\beta}
\newcommand{\ep}{\epsilon}  \newcommand{\vep}{\varepsilon}
  
    \newcommand{\la}{\lambda}
\newcommand{\sig}{\sigma}

\newcommand{\gam}{\gamma}
\newcommand{\R}{\mathbb{R}}
\newcommand{\N}{\mathbb{N}}

\newcommand{\ms}{\mathbb{S}}\newcommand{\ti}{\tilde }
\newcommand{\pt}{\partial_t}\newcommand{\pa}{\partial}
\newcommand{\les}{{\lesssim}}
\newcommand{\ges}{{\gtrsim}}

\newcommand{\hn}{\mathbb{H}^n}

\newcommand{\De}{\Delta}
\newcommand{\de}{\delta}

\newcommand{\hf}{\frac{1}{2}}

\newcommand{\one}{\uppercase\expandafter{\romannumeral1}}
\newcommand{\two}{\uppercase\expandafter{\romannumeral2}}
\newcommand{\three}{\uppercase\expandafter{\romannumeral3}}
\newcommand{\four}{\uppercase\expandafter{\romannumeral4}}
\newcommand{\five}{\uppercase\expandafter{\romannumeral5}}
\newcommand{\six}{\uppercase\expandafter{\romannumeral6}}
\newcommand{\sev}{\uppercase\expandafter{\romannumeral7}}
\newcommand{\eig}{\uppercase\expandafter{\romannumeral8}}

\newcommand{\pri}{\prime}
 
\newcommand{\ps}{\frac1{p}}

\title
[Critical John problem]
{Critical conditions of nonlinearities in the John problem}

\author{Chengbo Wang}
\address{School of Mathematical Sciences\\
Zhejiang Normal University\\
Jinhua 321004, P. R. China}
\email{wangcbo@zjnu.edu.cn}

\author{Xiaoran Zhang$^{*}$}\thanks{* Corresponding author}
\address{
Beijing Institute of Mathematical Sciences and Applications, Beijing 101408, P.R.China, 
Yau Mathematical Sciences Center, Tsinghua University, Beijing 100084, P.R.China}
\email{zhangxiaoran@bimsa.cn}

\thanks{
The first author was
supported by NSFC 12141102.
The second author was supported by
 CPSF 2025M783073. }

\begin{document}

\bibliographystyle{plain}

\begin{abstract}
In this paper, we determine the sharp criteria for the  nonlinearities in the John problem $\partial_{t}^2 u-\De_{\R^3}u= F(u)$ so that the problem admits global solutions for small initial data with compact support.
The criteria is of Dini type.
For the situation in which the criteria is not satisfied, we show that the solution will blow up in finite time for some small initial data, together with an almost sharp estimates of the lifespan.
For the proof, we mainly rely on the precise pointwise estimates of the solution, based on the fundamental solution, John's iteration argument as well as an adaption of the slicing method.
\end{abstract}



\maketitle


\section{Introduction}

In the seminal work \cite{John79}, F. John studied the problem of global existence vs blow up for the small amplitude Cauchy problem with (nontrivial) compactly supported, smooth data
\beeq\label{nlw}
\begin{cases}
\Box u
=F(u)\ ,\\
u(0,x)=u_0
,\ u_t(0,x)=u_1\ ,
\end{cases}
\eneq
where $u : I \times \R^3 \to \R$ is a real scalar field 
and $\Box u :=\pt^2 u-\Delta u$ is the d'Lambertian.
In  \cite[Corollary II]{John79}, it is shown that the problem with nontrivial data could not admit global solutions, when $F(u)\ge |u|^p$ for some $p\in (1,1+\sqrt{2})$, and $|F(u)|\le C|u|$ for $|u|\le 1$. 
On the other hand, in  \cite[Theorem II]{John79}, global existence for data with sufficiently small size is shown, whenever
$F(0)=F'(0)=F''(0)=0$ and $F\in C^{2,\theta}$ with some $\theta\in (\sqrt{2}-1, 1)$.
In particular, for the sample case of power functions, $F(u)=|u|^p$, the critical power $p_J=1+\sqrt{2}$ is determined, which
inspires the Strauss conjecture \cite{MR614228} for general spatial dimension $n\ge 2$. The conjectured critical power $p_S(n)$ is the positive root of 
$$(n-1)p^2-(n+1)p-2=0\ ,$$
with $p_S(3)=p_J=1+\sqrt{2}$. In the three dimensional case, it is also referred to as the John problem.

The Strauss conjecture has been systematically investigated. It is known that the problem admit global solution for any sufficiently small data, as long as $p\in (p_S(n), 1+4/(n-1)]$. On the other hand, such a result is not possible for $1<p\le p_S(n)$, in the sense that there exists a class of functions $(\phi,\psi)\in C_c^\infty(\R^n)$, so that the problem with data $(u_0,u_1)=(\ep \phi, \ep\psi)$ does not admit global solutions for any $\ep>0$.
  We refer interested readers to the classical works, Georgiev--Lindblad-Sogge  \cite{MR1481816}, Tataru  \cite{Ta01-2} for $p>p_S(n)$, Sideris \cite{MR744303} for $1<p<p_S(n)$, Schaeffer \cite{MR824205}, Yordanov--Zhang \cite{Yordanov}, Zhou \cite{MR2316656} for $p=p_S(n)$.
The classical result for the estimate of the lifespan $T(\ep)$ for $p=p_S(n)$, at least for $n\le 8$ or radial data, is
$$\log T(\ep)\sim \ep^{-p(p-1)}\,,$$
see the classical works Lindblad-Sogge \cite{lindblad-sogge1996}, Takamura-Wakasa\cite{takamura-wakasa2011},Zhou-Han\cite{zhou-han2014}.
For more history and recent works for the Strauss conjecture on asymptotically flat spacetimes, we refer the readers to a review by the first author \cite{Wang18}.

In the
field of nonlinear
  wave equations, it has been a general phenomenon that, the critical behavior could be likely extended beyond power functions, see Section \ref{sec:review} for a review on the history on various related model problems.
In this paper, we are interested in investigating the criteria for the general nonlinearity, beyond the 
scope of power functions, for the physical situation $n=3$.

\subsection{
Results for related model problems}\label{sec:review}
\subsubsection{Defocusing nonlinear wave equation}
One of the previous examples is the
 defocusing nonlinear wave equation $\Box u+|u|^{p-1}u=0$ in $\R^n$, the energy critical power is given by $1+4/(n-2)$, for which it is known that we have global well-posedness for $1<p\le 1+4/(n-2)$.
It was observed that 
the results could be extended for slightly supercritical nonlinearities,
in the form
$F(u)=|u|^{4/(n-2)}u\log^\al (2+u^2)$ or
$F(u)=|u|^{4/(n-2)}u\log^\al (\log (10+u^2))$ for some $\al>0$, see, e.g.,
 Tao \cite{MR2329385}, Roy \cite{MR2603799}, Shih \cite{MR3068539}, Bulut-Dodson \cite{MR4251268}.

\subsubsection{Damped wave equation}
Considering the damped wave equation
$\Box u+u_t=|u|^{p}$ in $\R^n$, similar to the Strauss conjecture,  Todorova--Yordanov
\cite{TY01}
found the critical power of the problem, given by $p_c=1+2/n$. 
Moreover, for the problem with 
the critical power $p_c=1+2/n$,
Zhang \cite{MR1847355} showed that the solutions also blow up
in a finite time.
See 
Lai--Zhou \cite{LaiZhou19} and references therein for more history, as well as the sharp lifespan estimates.

In \cite{Ebert},
Ebert, Girardi and Reissig tried to generalize the problem from the power nonlinearities to nonlinear function
$F(u)=|u|^{p_c}\mu(|u|)$, where 
$\mu: [0,\infty)\to [0,\infty)$ is a continuous, concave and increasing function satisfying $\mu(0)=0$.
More precisely, for the Cauchy problem of the damped wave equations
$$\Box u+u_t=|u|^{1+2/n}\mu(|u|)\ ,$$
small data global existence was shown for $n=1,2$, provided that
$\mu$ satisfies
$$\sup_{\la\in (0,1]}\sum_{1\le k\le n} \frac{\la^k |\mu^{(k)}(\la)|}{\mu(\la)}<\infty\ , $$
and the Dini condition
\beeq\label{dini-cond}
\int_0^{1}\frac{\mu(\la)}{\la}d\la <\infty\ .
\eneq
On the other hand,
suppose $F$ is convex and $\int_0^{1}\frac{\mu(\la)}{\la}d\la =\infty$, the solution blows up in finite time, whenever $u(0)=0$ and the  initial speed  is nontrivial and nonnegative.

\subsubsection{Strauss conjecture on hyperbolic spaces}
Another example is the analog of the Strauss conjecture on hyperbolic spaces, for which we consider the (shifted) wave equations on hyperbolic spaces $\mathbb{H}^n$,
$$\partial_{t}^2 u-(\De_{\hn}+(n-1)^2/4)u=|u|^p\ .$$
 In contrast with the Euclidean case, the geometry of $\mathbb{H}^n$ provides stronger decay and integrability at spatial infinity, in terms of the geodesic radius. 
Thanks to the improved dispersive estimates (see, e.g., Tataru \cite{Ta01-2}, Anker--Pierfelice--Vallarino \cite{APV12}), the Strauss conjecture on the hyperbolic space has been investigated in
Fontaine \cite{Fo97},
Sire--Sogge--Wang
 \cite{MR4026182}, where small data global existence is obtained for any sub-conformal powers $1<p<1+4/(n-1)$.
 In some sense, we can say that the critical power is given by $p_c=1$.
In other words, 
the critical behavior does not appear in the scale of the power-type nonlinearities. We remark also that similar global results are established for nontrapping (conformally compact) asymptotically hyperbolic manifolds, 
under certain natural spectral conditions, see Sire--Sogge--Zhang--Wang \cite{MR4169670,MR4432952}.

It turns out that the critical nature appears in the scale of logarithmic nonlinearities, at least for $\mathbb{H}^n$ with $n=2,3$, see
Wang--Zhang \cite{WZ25}
and Zhang
\cite{zhang2023critical}.
More precisely,
we investigate the critical behavior for the problem
$$
\partial_{t}^2 u-(\De_{\hn}+(n-1)^2/4)u=|u|(\sinh^{-1} 1/|u|)^{-(p-1)}\ ,$$
where the nonlinear function
$|u|(\sinh^{-1} 1/|u|)^{-(p-1)}$ behaves like $|u|(\log\frac{1}{|u|})^{-(p-1)}$  when $|u|\ll1$.
It is shown that we have
small data global existence for $p>3$, while blow up phenomenon appears for $p<3$.
That is,  the sharp critical exponent  is observed in the logarithmic scale of nonlinearity.

\subsubsection{Generalized Strauss conjecture}
As we have reviewed, considering \eqref{nlw} with $F(u)=|u|^p$ and $x\in \R^n$, the critical power is given by $p_S(n)$. In view of the related results for defocusing nonlinear wave equation, damped wave equation and the Strauss conjecture on hyperbolic spaces, it is a natural problem to investigate the similar problem with general nonlinearity, for which we refer as the generalized Strauss conjecture. 

The problem is initiated by Chen and Reissig in
\cite{Chen_2024}, where the authors
 considered
  $$\Box u=|u|^{p_S(n)}\mu(u)\ ,$$
where, similar to   \cite{Ebert},
$\mu: [0,\infty)\to [0,\infty)$ is a continuous, concave and increasing function satisfying $\mu(0)=0$.

It was  shown that if $n\ge 2$ and 
\beeq\label{C-R1}\lim_{s\to 0+} \mu(s)(\log 1/s)^{1/p_S(n)}\gg 1\ ,\eneq
the solution blows up in finite time for the problem with nonnegative, nontrivial and compactly
supported smooth data.
To indicate the (possibly) sharpness of the condition \eqref{C-R1}, for spatial dimension $n=3$, global existence with small radial data was proven, provided that
\beeq\label{C-R2}\lim_{s\to 0+} \mu(s)(\log 1/s)^{1/p_S(3)}=0\ ,\eneq
and
\beeq\label{C-R3} \mu(s)(\log 1/s)^{1/p_S(3)}(\log\log 1/s)\les 1\ ,\ \forall s\ll 1\ .\eneq

In particular, in three space dimensions, for the problem with sample logarithmic nonlinearity
\beeq\label{C-R}
\Box u=|u|^{p_S(3)}\left(\log\frac1{|u|}\right)^{-\gam}
=|u|^{1+\sqrt{2}}\left(\log\frac1{|u|}\right)^{-\gam}
\ ,
\eneq
 the critical exponent is determined to be $\gam_c=1/p_S(3)=\sqrt{2}-1$.

\subsubsection{Generalized Glassey conjecture}
Let $n\ge 2$ and $p>1$.
Consider the wave equation with nonlinearity depending on the derivative of the unknown: $$\Box u=|\pt u|^p\ ,$$
it is known that, for any nontrivial compactly supported smooth data,
the solution blows up in finite time, whenever $(n-1)(p-1)\le 2$.
On the other hand, for $p>1+2/(n-1)$, we have small data global existence when $n\le 3$. The high dimensional case is known for
$p\in (1+2/(n-1), 1+2/(n-2))$ and radial data. See
Zhou \cite{Zh01},
Hidano--Wang--Yokoyama \cite{HWY2} and references therein.


With the critical power $p_c=1+2/(n-1)$ and
$F(s)=|s|^{p_c}\mu(|s|)$, the generalized Glassey conjecture is recently studied, by
Chen and Palmieri \cite{chenGlassey} and
Shao \cite{Shao24}. More precisely, let $\mu: [0,\infty)\to [0,\infty)$
be a continuous,
increasing, concave function and satisfies $\mu(0) = 0$, the nonlinear wave equations under consideration is 
$$\Box u=F(\pt u)=|\pt u|^{1+2/(n-1)}\mu(|\pt u|)\ .$$
In \cite{chenGlassey}, finite in time blow up is shown when
$F$ is convex and the Dini condition
\eqref{dini-cond} is not satisfied.
On the other hand, when $n=3$ and the initial data are radial, small data global existence is proved when $\mu$ satisfies the Dini condition
\eqref{dini-cond} and
$$\sup_{\la\in (0,1]}\sum_{1\le k\le 2} \frac{\la^k |\mu^{(k)}(\la)|}{\mu(\la)}<\infty\ .$$
The nonradial global results for $n=2, 3$, as well as the analogs for the problem in nontrapping exterior domains, were obtained later in Shao \cite{Shao24}, by the local energy estimates and Klainerman-Sobolev inequalities. 

\subsection{Investigate the sharp critical condition}
After reviewing the history on various related model problems, we find that, at least for the damped wave equations or the generalized Glassey conjecture, the sharp condition for the problem to admit small data global solutions has been obtained for the sample spatial dimensions, which is basically given by the Dini condition \eqref{dini-cond}.

However, for the generalized Strauss conjecture in spatial dimension three, the sharp condition is only known for the  logarithmic nonlinearities
\eqref{C-R} with critical power $\gam_c=1/p_S(3)=1/p_J$, which generalizes John's result \cite{John79}.

Concerning the sharp criteria for the nonlinear functions, it remains unclear whether it is given by the logarithmic type nonlinearities or Dini type conditions.
This suggests that the generalized Strauss conjecture is much more delicate.

In this paper, we try to determine the criteria for the  generalized Strauss conjecture in  the physical dimension $n=3$, which is corresponding to the John problem.
In this subsection, we would like to illustrate what is the natural expected criteria, based on the fundamental solution and John's iteration argument.


%
 To determine the  critical condition of $\mu$ for \eqref{nlw} with
\beeq\label{nonlin}F(u)=|u|^{1+\sqrt{2}}\mu(|u|)\ ,\eneq
 we start from pointwise estimates of the radial solution in $\R^3$. It is well known that the linear radial solution admits a particularly simple expression in $\mathbb{R}^3$. By Duhamel's principle, we can rewrite the integral formulation for the solution to \eqref{nlw}, with compactly supported radial data:
\beq
\begin{split}
u(t,r)&= u^0+LF(u)=  u^0(t,r)+\frac1{2r}\int_0^t\int_{|t-s-r|}^{t-s+r}(|u|^{1+\sqrt{2}}\mu(|u|))(s,\la)\la d\la ds\\
&= u^0(t,r)+\frac1{4r}\int_{|t-r|}^{t+r}\int_{-\be}^{t-r}(|u|^{1+\sqrt{2}}\mu(|u|))(\al,\be)\frac{\be-\al}2d\al d\be\ ,
\end{split}
\eeq
where $u^0$ is the homogeneous solution, and $\al=s-\la$, $\be=s+\la$ are the double null coordinates.

By sharp Huygens' principle and the bootstrap argument, global existence could be proved, provided that 
 the following pointwise estimates could be kept under iteration, 
\beeq\label{intr-u0}|u(t,r)|\les 
\<t+r\>^{-1}h(t-r)\ ,\eneq
for some function $h$ to be determined. Here
we denote $\<x\>:=2+|x|$.
It turns out that we may choose
\beeq\label{intr-u}|u(t,r)|\les 
\<t+r\>^{-1}\<t-r\>^{1-\sqrt{2}}\mu(\<t-r\>^{-1})\les \<t+r\>^{-1}\ .\eneq
Then by the iteration $u^n= u^0+LF(u^{n-1})$, we have
\beq
|LF(u)(t,r)|\les \frac{1}{r}\int_{|t-r|}^{t+r}\<\be\>^{-\sqrt{2}}\mu(\<\be\>^{-1})\int_{-\be}^{t-r}\<\al\>^{-1}\mu^{1+\sqrt{2}}(\<\al\>^{-1})d\al d\be\ .
\eeq
For the purpose of closing the iteration, 
 this appears to be closely related to the
 integrability of the integral
$$ \int_{-\be}^{t-r}\<\al\>^{-1}\mu^{1+\sqrt{2}}(\<\al\>^{-1})d\al\ ,$$
which is reduced to the  following condition 
\beeq\label{global-mu}\int_0^{1}\frac{\mu^{1+\sqrt{2}}(\la)}{\la}d\la<\infty\ .\eneq
We observe that this is a Dini type condition which generalizes the logarithmic form in \eqref{C-R2}. This suggest it might be a good candidate for the sharp condition.

On the other hand, to illustrate the blowup phenomenon, 
we may want to see increasing lower bound for the solution, by the John's iteration argument. Beginning with certain well-chosen  radial, compactly supported, smooth initial data, we could have
$$u^0(t,r)\ges \<t+r\>^{-1}$$
near the light cone $0< t-r\le 1$. Thanks to the nonnegative nature of the fundamental solution and nonnegative nonlinearity, we have 
$u\ge u^0(t,r)$. Accordingly, applying the first iteration yields
$$u(t,r)\ges \frac1{r}\int_{t-r}^{t+r}\be^{-\sqrt{2}}\mu(\be^{-1})d\be\ges(t+r)^{-1}(t-r)^{1-\sqrt{2}}\mu((t-r)^{-1})\ges (t+r)^{-2}\ ,$$
for $t-r\ge 2$, see Proposition \ref{pointwise-low}. Then the second iteration arrives at
$$u(t,r)\ges \frac1{2r}\int_{t-r}^{t+r}\be^{-1-\sqrt{2}}(\be-(t-r))\mu(\be^{-2})d\be\int_{2}^{t-r}\al^{-1}\mu^{1+\sqrt{2}}(\al^{-1})d\al\ .$$
In the sample case where $\mu(\be^{-2})\sim \mu(\be^{-1})$, we see that the last term will surpass the original lower bound, provided that
$$\int_{2}^{t-r}\al^{-1}\mu^{1+\sqrt{2}}(\al^{-1})d\al\gg 1\ .$$
This is exactly the case when $\mu$ is violating \eqref{global-mu}.
It suggests that blow up might also closely related with this critical integral.

Based on the preceding analysis, we infer that 
\eqref{global-mu} is the sharp condition for the problem \eqref{nlw} to admit small data global solutions,
under certain natural restrictions on $\mu$.

\subsection{Our results}
In this manuscript, we prove that  \eqref{global-mu} is exactly the sharp condition for the problem \eqref{nlw} to have small data global existence.
For simplicity of the presentation, we shall simply write $p=1+\sqrt{2}$ as the John-Strauss critical power hereinafter.

Let us first state our positive result.
\begin{thm}[Global existence]\label{global-thm}
Let $\mu: [0,\infty)\to [0,\infty)$ be a continuous and increasing function such that
$F(u)=|u|^p\mu(|u|)\in C^2$ and the Dini type condition  \eqref{global-mu} holds.
Then the Cauchy problem \eqref{nlw} enjoys small data global existence. 
More precisely,
 for any $(f,g)\in C^3_c\times C^2_c$, there exists $\vep_0>0$ so that the problem  \eqref{nlw} with data $(u_0,u_1)=(\vep f, \vep g)$ admits a global classical solution $u\in C^2(\R\times\R^3)$ for any $|\vep|<\vep_0$.
\end{thm}

On the other hand, to illustrate the sharpness of 
the Dini type condition  \eqref{global-mu}, 
we consider the classical solution to the Cauchy problem
\beeq\label{blow-eq}
\Box u=F(u):=|u|^{p}\mu(|u|)\ ,\ u(0,x)= 0,\ u_t(0,x)= u_1\ , x\in\R^3 .\eneq

\begin{thm}[Blow up and estimates of the lifespan]\label{blowup-thm}
Let $\mu: [0,\infty)\to [0,\infty)$ be a continuous and increasing function that violates the Dini type condition  \eqref{global-mu}, that is,
\beeq\label{blowup-mu}
\int_0^{1}\frac{\mu^{p_S}(\la)}{\la}d\la=\infty\ .
\eneq
Suppose also that $F\in C^2$,
 $\mu\in C^1((0,\infty); (0,\infty))$ and
\beeq\label{diff-poly}
\frac{\la \mu^{\pri}(\la)}{\mu(\la)}=o(1)\ ,\ as\ \la\to0_+.
\eneq
Then
 the classical solution $u$ to the Cauchy problem \eqref{blow-eq} blows up in finite time, for any radial, non-negative and non-trivial data $u_1\in C_c^2$. More precisely, for any radial, non-negative and non-trivial $C_c^2$ function $f$, and any $0<\de\ll1$, there exist $\vep_0=\vep_0(f,\de)$ and positive constants $c,C$ independent of $\vep,\de$, such that for $\vep\in(0,\vep_0]$,
 the lifespan $T(\vep)$ is estimated by
 \beeq\label{eq-life}c\vep^{-p(p-1)}\le A_{\vep}(T(\vep))\le C\vep^{-p(p-1)-\de/p}\eneq
 with $u_1=\vep f$, where $A_{\vep}$ is defined in \eqref{def-A} in Section \ref{pre}.

\end{thm}
\begin{remark}
    Notice that the Dini-type condition \eqref{global-mu}, together with the continuity of $\mu$ in Theorem \ref{global-thm}, already implies $\mu(0)=0$ in this case. On the other hand, if \eqref{global-mu} fails, as in Theorem \ref{blowup-thm}, then $\mu(\la)>0$ for any $\la>0$.
\end{remark}

\begin{remark}
  Note that $\mu\equiv1$ is also a typical example for Theorem \ref{blowup-thm}, and there is no such loss $\de$ in the estimate of the lifespan, which is consistent with the
   well known critical upper bound, see
   Takamura-Wakasa\cite{takamura-wakasa2011}, Zhou-Han\cite{zhou-han2014}.
  See Remark \ref{mu=1} for details.
\end{remark}

\begin{remark}\label{ex-mu}
Compared with the results established by Chen and Reissig \cite[Theorem 2.1, 2.2]{Chen_2024}, our work establishes the sharp critical condition  \eqref{global-mu}.
Specifically, we would like to give typical examples of $\mu$ in Theorem \ref{blowup-thm} for $\la\ll1$:
\beeq\label{eg.2} \mu_k(\la)= \left(\log\frac1\la\right)^{-\frac1{p}}\left(\log\log\frac1\la\right)^{-\frac1{p}}\cdots \left(\log^k\frac1\la\right)^{-\frac1{p}}\ ,
\eneq
Here, we use $\log^k$ to represent $k$-th composition of $logarithmic$ mapping.
In particular, Theorem \ref{blowup-thm}  with $\mu_1$ gives us the blow up result for 
\eqref{C-R} with the critical power.
Remark also that, for any $k\ge 1$, $\mu_k$ satisfies an additional property $\mu(\la)\sim\mu(\la^2)$. For such functions, we claim that
the proof of Theorem \ref{blowup-thm}  will be much easier, and we refer the readers to a previous version of this manuscript to the proof for $\mu_k$, see arXiv:2405.12761 v1.
   However, such a property is not enjoyed for general $\mu$
  satisfying \eqref{diff-poly} while 
   violating  \eqref{global-mu}. The absence of this property poses one of the major difficulties in the analysis. To handle this issue, we introduce a key function $A$, which is naturally connected to the integral
   \eqref{global-mu}, see Section \ref{pre}.
\end{remark}

Concerning the idea of the proof, 
 we mainly rely on the precise pointwise estimates of the solution, based on the fundamental solution and John's iteration argument. 
 The proof of the global well-posedness result in Theorem \ref{global-thm} is presented in Section \ref{proof-global}.
The proof of the blow up part in Section \ref{blow-up-upper-life} is more delicate.
The functional method is usually used in the blow up analysis, instead, we turn to a precise pointwise estimate here, in a fashion similar to that of John \cite{John79}. In this way, we can make the blow up nature more clear.
More precisely, by exploiting the positive fundamental solution, 
we obtain an improved initial bound  \eqref{lower-2} inside the light cone in Section \ref{initial-lo-bd}.
 To facilliate the iteration, we introduce a key function $A$ in Section \ref{pre} -\eqref{def-A}, which is naturally connected to the integral
   \eqref{global-mu}.
    Based on 
an adaption of
the slicing method of Agemi-Kurokawa-Takamura \cite{MR1785116}, we obtain a sequence of lower bounds \eqref{iter-lowerbdform} for the solution, which is increasing and diverges to infinity in finite time in Section \ref{iter-proc}. It is also precise enough for us to obtain an upper bound estimate of the lifespan. Finally, in Section \ref{lower-esti-life}, we complete the almost sharp estimate by establishing the lower bound of the lifespan.

Before concluding the introduction, let us take an outlook to the general spatial dimension.
As we have recalled, the blow up result is also known \cite{Chen_2024} for general $n\ge 2$, under the condition
\eqref{C-R1}. In view of the sharpness of the Dini type condition \eqref{global-mu} for $n=3$, as well as the sharpness of the Dini condition  \eqref{dini-cond} for
the damped wave equations and the generalized Glassey conjecture, it seems natural to expect that it is also the sharp condition for general $n$.

In conclusion, we
propose the following conjecture:
\begin{conjecture}\label{our-conj}
Consider the Cauchy problem \eqref{nlw} with general spatial dimension $n\ge 2$ and nonlinearity 
$$F(u)=|u|^{p_S(n)}\mu(|u|)\ ,$$
 where $\mu: [0,\infty)\to [0,\infty)$ is a continuous and increasing function. Suppose in addition that 
 $\mu\in C^k((0,\infty); (0,\infty))$ for some $k\ge 1$ and
\beeq\label{diff-poly2}
\frac{\la^j \mu^{(j)}(\la)}{\mu(\la)}=o(1)\ ,\ as\ \la\to0_+\ , \forall 1\le j\le k\ .
\eneq
Then we have  small data global existence if and only if the following Dini type condition is satisfied:
\beeq\label{eq-criteria}
\int_0^{1}\frac{\mu^{p_S}(\la)}{\la}d\la<\infty.
\eneq
\end{conjecture}
For the high dimensional case, we may need to suppose some additional conditions on $\mu(\la)$ for large $\la$. However, for $n=2$, we believe the assumption is sufficient for the conjecture.

\subsection*{Notation} 
\begin{itemize}
\item We use $A\lesssim B$ to denote
$A\leq CB$ for some fixed large positive constant C which may vary from line to
line, and similarly we use $A\ll B$
to denote $A\leq C^{-1} B$. We employ $A\sim B$ when $A\lesssim
B\lesssim A$.  Particularly, when $C=C_k$ only depends on parameter $k$, we use $A\les_k B$ and $A\sim_k B$.
\item We denote $\<x\>:=2+|x|$.
\item We denote $x\vee y:=\max\{x,y\}$ and $x\wedge y:=\min\{x,y\}$.
\end{itemize}

\section{Preliminary}\label{pre}

Considering the linear wave equations \beeq\label{wave}
\begin{cases}
\partial_{t}^2 u-\De u=F\ , x\in\R^3\ ,\\
u(0,x)= u_0, u_t(0,x)= u_1\ ,
\end{cases}
\eneq
Duhamel's principle tells us that
 \eqref{wave} is equivalent to the integral equation
\beeq\label{ie}
u(t,x)= u^0(t,x)+(LF)(t,x)= u^0(t,x)+\int_0^t S(\tau)F(t-\tau, x)d\tau,
\eneq
where $u^0= \pt S(t)u_0(x)+  S(t)u_1(x)$, and $S(t)u_1(x)$ is the solution for the linear homogeneous equation with data $u_0=0$.

Based on the explicit fundamental solution $\frac{1}{2\pi}\de(t^2-|x|^2)$,
 we can easily obtain the following classical dispersive estimate for the homogeneous solution $u^0$.


\begin{lem}[Linear estimates]\label{linear-esti}
Let $n=3$  and $r=|x|$. Suppose the initial data $(u_0, u_1)\in C^3_c\times C^2_c$ are supported in $B_R(O)$. Then
there exists $N>0$ so that 
\beeq\label{eq-linear}
|u^0(t,x)|\le  N\<t+r\>^{-1}\chi_{|t-r|\le R}\ ,
\eneq
for the solution to  \eqref{wave} with $F=0$.
\end{lem}


Next we establish the preliminary of the proof of Theorem \ref{blowup-thm}.
Firstly, based upon the assumption made on $\mu$ in Theorem \ref{blowup-thm}, we derive a key property of this function.
\begin{lem}\label{lem-la}
Let $\mu\in C^1((0,\infty); (0,\infty))$ be an increasing function  satisfying 
\eqref{diff-poly}, that is,
$$\lim_{\la\to 0+}\frac{\la \mu^{\pri}(\la)}{\mu(\la)}=0\ .$$
Then
\begin{itemize}
\item[(1)]
 for any pair $(a,b)\in \R^2$ with $a\neq 0$, 
\beeq\label{deri-lamu}
\lim_{\la\to 0+}\frac{\frac{d}{d\la}(\la^a\mu^b(\la))}{a \la^{a-1}\mu^b(\la)}=1
\eneq
\item[(2)]
For any $\de_0>0$, there exists $\la_0\in (0,1)$, so that
\beeq\label{lowerbd-mu}
\mu(\theta \la)\le\mu(\la)\le \theta^{-\de_0}\mu(\theta\la),\ \forall \la\in(0,\la_0], \forall  \theta\in (0,1].
\eneq
In particular, 
we have 
\beeq\label{lowerbd-mu2}
\mu(\la)\ge \la^{\de_0}\la_0^{-\de_0}\mu(\la_0),\ \forall \la\in (0, \la_0]\ .\eneq
\end{itemize}

\end{lem}


\begin{proof}
The calculation for \eqref{deri-lamu} is straightforward,
as we have
$$\frac{\frac{d}{d\la}(\la^a\mu^b(\la))}{a \la^{a-1}\mu^b(\la)}=1+\frac{b}a \frac{\la \mu'(\la)}{\mu(\la)}\ .$$
As for \eqref{lowerbd-mu}, we make a change of variables.
Let $\theta=e^{-\de}$, and
 $f(\la)=\frac{\mu(\la)}{\mu(\theta\la)}$, we have
$\log f(e^s)=\log \mu(e^s)-\log \mu(e^{s-\delta})$.
Observing that $$\frac{d}{ds}\log \mu(e^s)=\frac{e^s \mu'(e^s)}{\mu(e^s)}\ ,$$
then for any $\de_0>0$, there exists $\la_0\in (0,1)$ so that for any $\la=e^s<\la_0$,
we have $\frac{d}{ds}\log \mu(e^s)\in [0, \de_0]$. Thus, integrating from $s-\delta$ to $s$ gives us
$$
0\le \log f(e^s)\le \de_0 \de=-\de_0\log \theta\ , \forall s<\log \la_0\ .
$$ Returning to the original variable, we get
$$1\le f(\la)\le \theta^{-\de_0}\ ,\  \forall \la<\la_0\ ,$$
which gives us the desired inequality  \eqref{lowerbd-mu}.
\end{proof}

Based on Lemma \ref{lem-la} $(2)$, we obtain the equivalence  $$\mu(\la)\sim_C\mu(C\la)$$ for small $\la$. However,
we emphasize that $\mu$ does not satisfy the stronger equivalence
\beeq\label{fail-mu}
\mu(\la)\sim\mu(\la^2),\quad \la\ll1,\eneq
which was  crucial for the blow-up iteration in the previous works.

In any case, to connect between the scales $\la$ and $\la^2$, we introduce the following function, which is closely related to $\mu$ and the growth induced by the blow-up assumption \eqref{blowup-mu}.
For given $\de>0$, we denote
\beeq\label{def-A}
A_\de(t)=\int_{t^{-1}}^{1}\frac{\mu^{p}(\de\la)}{\la}d\la.
\eneq
Clearly, $A_\de:[1,+\infty)\to[0,+\infty)$ is an increasing function, $A_\de(1)=0$ and tends to infinity as $t\to\infty$. We observe some critical properties of $A_\de$ as follows.

\begin{lem}\label{lem-A}
Let 
$\mu\in C ([0,\infty); [0,\infty))$
 be an increasing function, then we have the following results:
\begin{itemize}
\item[(1)] $A_\de^{\pri}(t)=\frac{\mu^{p}(\de/t)}t$.
\item[(2)]  For any $\theta\in (0,1]$, there exists some positive constant $C_{\theta}$ such that for $t>2/\theta$,
$$A_\de(\theta t)\le A_\de(t)\le C_{\theta}A_\de(\theta t)\,.$$
\item[(3)] For $t>1$, $$A_\de(t)\le A_\de(t^2)\le 2 A_\de(t)\,.$$
\end{itemize}
\end{lem}

\begin{proof}
The first identity is trivial. Then we only consider the remaining two inequalities.
Firstly, a brief calculation shows that
\beq
\begin{split}
A_\de(\theta t)\le A_\de( t)&=\int_{1/t}^{1}\frac{\mu^{p}(\de\la)}{\la}d\la=A_\de(\theta t)+\int^{1/(\theta t)}_{1/t}\frac{\mu^{p}(\de\la)}{\la}d\la\\
&\le A_\de(\theta t)+ (\log 1/\theta)\mu^{p}(\de/(\theta t))\ .
\end{split}
\eeq
Recall
that, if $t>2/\theta$, $A_\de(\theta t)\ge\int_{1/(\theta t)}^{2/(\theta t)}\frac{\mu^p(\de \la)}{\la}d\la \ge (\log 2) \mu^p(\de/(\theta t))$, we get
$$A_\de(\theta t)\le A_\de( t)\le \left(1-\frac{\log \theta}{\log 2}\right)A_\de(\theta t)\ ,$$
which is the second inequality.

For the third one, we make a change of variable in the integral, with $\la=s^2<1$, and using the fact that
$\mu^p(\de s^2)\le \mu^p(\de s)$, we obtain
\begin{align*}
A_\de(t)\le A_\de(t^2)&=\int_{t^{-2}}^1\frac{\mu^{p}(\de \la)}{\la}d\la\\
&=2\int_{t^{-1}}^1\frac{\mu^{p}(\de s^2)}{s}d s\\
&\le 2\int_{t^{-1}}^1\frac{\mu^{p}(\de s)}{s}d s=2 A_\de(t)\ ,
\end{align*}
which completes the proof. 

\end{proof}

\section{Global existence: proof of Theorem \ref{global-thm}}\label{proof-global}
In this section, we give the proof of Theorem \ref{global-thm}. 
For sufficiently small $\vep>0$,
by
 \cite[Chapter \one, Theorem 5.1]{So08},
for some $T>0$,  we have the classical local solution $u\in C^2([0,T]\times\R^3)$ due to $F\in C^2$ and $F(0)=0$. Moreover, the corresponding continuation criterion shows that a finite maximal lifespan can occur only if the solution becomes unbounded. Therefore, to obtain global existence, it suffices to prove that there exists a constant $C$ such that
$$\sup_{x}|u(t,x)|\le C\ ,$$
for all $t$ in the lifespan of the solution.

Since $\mu$ is continuous on $[0,\infty)$, it is bounded on $[0,1]$. Absorbing this bound into the implicit constants, we may assume without loss of generality that
$$\mu(s)\le1,\quad 0\le s\le1,$$
and hence $\mu(\<y\>^{-1})\le1$. 

By Lemma \ref{linear-esti} and without loss of generality we assume $R=1$, there exists a constant $C_0>0$ such that
$$|u^0(t,x)|\le \vep C_0\<t+r\>^{-1}\<t-r\>^{-\ps}\mu(\<t-r\>^{-1})=: \vep C_0\Phi_0(t,r)\ .$$

Then bootstrap argument tells us that the proof is reduced to establishing the following boundness estimate.
\begin{lem}\label{bound-lem}
Under the assumption of Theorem \ref{global-thm}, there exists $\vep_0>0$ so that for any $\vep\in(0,\vep_0]$, we have 
\beeq\label{bound}
|LF(u)(t,x)|\le \frac{\vep}{2}  C_0\Phi_0(t,r)\ ,
\eneq
whenever $|u(t,x)|\le 2\vep C_0\Phi_0(t,r)$.
\end{lem}
\begin{proof}

By Duhamel’s principle and the well-known solution formula of linear wave equation for radial data
$$u^0(t,x)=S(t)u_1+\pt S(t)u_0\ ,$$
\beeq\label{radial-fom}
S(t)u_1(x)=tM^tu_1(x)=\frac{t}{|\ms^2|}\int_{\ms^2}u_1(x+ty)d\sig_y=\frac1{2r}\int_{|t-r|}^{t+r}u_1(\la) \la d\la\ ,
\eneq
where $|\ms^2|$ denotes the area of the sphere $\ms^2$. Thus, \eqref{ie} and \eqref{radial-fom} imply 
\beq
\begin{split}
|LF(u)(t,x)|&\le\frac1{2r}\int_0^t\int_{|t-s-r|}^{t-s+r}F(2\vep C_0\Phi_0)(s,\la)\la d\la ds\\
&\le\frac{(2\vep C_0)^{p}}{2r}\int_0^t\int_{|t-s-r|}^{t-s+r}\<s+\la\>^{-p}\<s-\la\>^{-1}\mu^{p}(\<s-\la\>^{-1})\mu(\Phi_0)\la d\la ds,
\end{split}
\eeq
when   $2\vep C_0<1$.
Introducing new variables of integration 
\beeq\label{new-vari}
\al=s-\la, \be=s+\la, 
\eneq
we have
\beq
|LF(u)(t,x)|\le\frac{(2\vep C_0)^{p}}{8r}\int_{|t-r|}^{t+r}\int_{-\be}^{t-r}\<\be\>^{-p}\<\al\>^{-1}\mu^{p}(\<\al\>^{-1})\mu(\Phi_0(\al,\be))(\be-\al)d\al d\be.
\eeq
Since $\Phi_0(\al,\be)\le \<\be\>^{-1}\le\<t-r\>^{-1}$, 
there exists a constant $C_1>0$ such that
\beq
|LF(u)(t,x)|\le\frac{C_1\vep^{p}}{r}\mu(\<t-r\>^{-1})\int_{|t-r|}^{t+r}\<\be\>^{1-p}\int_{-\be}^{t-r}\<\al\>^{-1}\mu^{p}(\<\al\>^{-1})d\al d\be\ .
\eeq

Under the assumption \eqref{global-mu} of $\mu$, we can control the dominant term:
$$
\int_{-\be}^{t-r}\<\al\>^{-1}\mu^{p}(\<\al\>^{-1})d\al
\les\int_{0}^{\be}\<\al\>^{-1}\mu^{p}(\<\al\>^{-1})d\al
\les\int_0^{1}\frac{\mu^{p}(\la)}{\la}d\la<\infty\ .
$$
Then we get for some $C_2>0$,
\beq
|LF(u)(t,x)|\le C_2\vep^{p}\mu(\<t-r\>^{-1})\frac1r\int_{|t-r|}^{t+r}\<\be\>^{1-p}d\be\ .
\eeq
When $r\le t/2$, that is, $t+r\le 3(t-r)$,the integral can be handled by exploiting the monotonicity of the integrand:
$$\frac1r\int_{t-r}^{t+r}\<\be\>^{1-p}d\be\le 2\<t-r\>^{1-p}\le 6\<t+r\>^{-1}\<t-r\>^{2-p}\ .$$
On the other hand, 
if $t\le 2$, we have $r\le t+1\le 3$,
$\<t+r\>\le 7\le 4\<t-r\>$,
similar argument gives us
$$\frac1r\int_{|t-r|}^{t+r}\<\be\>^{1-p}d\be\le \frac{2(r\wedge t)}{r}\<t-r\>^{1-p}\le 8\<t+r\>^{-1}\<t-r\>^{2-p}\ .$$
For the remaining case
 $r\ge t/2\ge 1$, we have $ \<t+r\>\le   3\<r\>\le 9 r$ it follows directly from integration by parts:
$$\frac1r\int_{|t-r|}^{t+r}\<\be\>^{1-p}d\be
=\frac{\<t-r\>^{2-p}-\<t+r\>^{2-p}}{(p-2)r}
\le \frac{9}{p-2}\<t+r\>^{-1}\<t-r\>^{2-p}\ .$$
With the fact of $2-p=-\ps$, we conclude the proof by taking $\vep_0$ to be sufficiently small.
\end{proof}

\section{Blow up and upper bound of the lifespan}\label{blow-up-upper-life}

Under the assumption of initial data $u_1=\vep f\in C_c^2$ and
 $F\in C^2$, 
 for sufficiently small $\vep>0$, there exists $T>2$ so that
we have a unique classical local solution $u\in C^2([0,T]\times\R^3)$.
 The supremum of all such $T$ is called the lifespan of the solution,
  denoted by $T(\vep)$.
  Without loss of generality, we will assume that the support of $u_1=\vep f$ is in the unit ball $B_1(O)$.

In this section, we give estimates of the lifespan in an implicit 
expression, where the upper bound of the lifespan implies that 
the solution blows up in finite time for any small $\vep>0$. 


\subsection{Initial lower bound}\label{initial-lo-bd}
From the expression of solutions \eqref{ie} and \eqref{radial-fom}, 
we deduce the following pointwise lower bound estimate. 

\begin{prop}\label{pointwise-low}
Under the assumption of Theorem \ref{blowup-thm}, there exist constants $r_0, C, C_0\in(0,1)$ depending only on $f$ such that for sufficiently small $\vep>0$, the classical solution of the Cauchy problem \eqref{blow-eq} has the pointwise lower bound 
\beeq\label{lower-1}
u(t,r)\ge \frac{C \vep}{t+r},\  \forall 0<t-r<r_0<1<t+r\ ,  
\eneq
and
\beeq\label{lower-2}
u(t,r)\ge C_0\vep^p(t+r)^{-1}(t-r)^{-\ps}\mu\left(\frac{\vep}{t-r}\right), \ \forall t-r>2\ .
\eneq
\end{prop}
\begin{proof}
  As $u_1$ is supported in the unit ball $B_1(O)$,
by \eqref{ie} and \eqref{radial-fom}, 
we have the following lower bound for $t>r$
$$u(t,r)\ge u^0(t,r)= \frac{\vep}{2r}\int_{t-r}^{t+r}f(\la)\la d\la\ .$$
Since $f$ is non-trivial and non-negative, there exists $r_0\in(0,1)$ such that
$$\int_{r_0}^{1}f(\la)\la d\la >0\ .$$
Then, for $0<t-r<r_0<1<t+r$, we have
\beeq\label{initial-bd}
u(t,r)\ge \frac{\vep}{t+r} \int_{r_0}^{1}f(\la)\la d\la\ge \frac{C\vep}{t+r}\ ,
\eneq
for any $C\in (0,1)$ with $C\le \int_{r_0}^{1}f(\la)\la d\la$.

Subsequently, we employ the following iterative form with the new coordinate $(\al,\be)$ in \eqref{new-vari}:
\beeq\label{iter-lb}
\begin{split}
u(t,r)&\ge\int_0^tS(t-s)F(u)(s)ds=\frac1{2r}\int_0^t\int_{|t-s-r|}^{t-s+r}(|u|^{p}\mu(|u|))(s,\la)\la d\la ds\\
&=\frac1{8r}\int_{t-r}^{t+r}\int_{-\be}^{t-r}(|u|^{p}\mu(|u|))(\al,\be)(\be-\al)d\al d\be.
\end{split}
\eneq
Then the initial lower bound \eqref{initial-bd} yields 
\beq
\begin{split}
u(t,r)&\ge \frac{C^p\vep^p }{8r}\int_{t-r}^{t+r}\int_0^{r_0} \be^{-p}\mu(C\vep\be^{-1})(\be-\al)d\al d\be\\
&\ge \frac{C^p\vep^p r_0}{16r}\int_{t-r}^{t+r}\be^{1-p}\mu(C\vep\be^{-1})d\be\ ,
\end{split}
\eeq
for $t-r>2$,
where in the last inequality we used the fact that $\be-\al\ge \be/2$ for $\al\in[0,r_0]$ and $\be>t-r>2>2r_0$.

When $t-r>2$ and $r<t/2$, which means $t-r\le t+r\le 3(t-r)$, thanks to Lemma \ref{lem-la}, it is clear that for sufficiently small $\vep$,
$\be^{1-p}\mu(C\vep\be^{-1})$ is a decreasing function by \eqref{deri-lamu} with $a=p-1,b=1$, and $$\mu(C\vep (t+r)^{-1})
\ge\mu(\frac C 3 \vep (t-r)^{-1})
\ge \frac{C}3 \mu(\vep (t-r)^{-1})$$ by \eqref{lowerbd-mu}. Thus, we have
\begin{eqnarray*}
u(t,r) & \ge& \frac{r_0C^p}8 \vep^p(t+r)^{1-p}\mu(C\vep(t+r)^{-1})\\
 & =& \frac{r_0C^p}8 \vep^p
 (t+r)^{-1}(t+r)^{2-p}\mu(C\vep(t+r)^{-1})\\
&\ge& \frac{r_0C^{p+1}}8 3^{1-p}
\vep^p(t+r)^{-1}(t-r)^{2-p}\mu(\vep(t-r)^{-1})\ .
\end{eqnarray*}

On the other hand,  by \eqref{deri-lamu} with $a=p-2,b=1$,  we get
$$\be^{1-p}\mu(C\vep\be^{-1})\ge
\frac 12 
\pa_{\be}\frac{\be^{2-p}\mu(C\vep\be^{-1})}{2-p}$$
for sufficiently small $\vep$. Hence,
for $t/2\le r<t$, we have
$t+r\ge 3(t-r)$, and thus
\beq
\begin{split}
u(t,r)&\ge \frac{r_0C^p\vep^p}{32r}\int_{t-r}^{t+r}\pa_\be\left(\frac{\be^{2-p}\mu(C\vep\be^{-1})}{2-p}\right)d\be\\
&\ge \frac{r_0C^p\vep^p}{32(p-2) t}\left[(t-r)^{2-p}\mu(C\vep(t-r)^{-1})-(t+r)^{2-p}\mu(C\vep(t+r)^{-1})\right]\\
&\ge \frac{r_0C^p\vep^p}{32(p-2) t}\left[(t-r)^{2-p}-(t+r)^{2-p}\right]\mu(C\vep(t-r)^{-1})\\
&\ge\frac{r_0C^{p+1}(1-3^{2-p}) }{32(p-2)}  \vep^p(t+r)^{-1}(t-r)^{2-p}\mu(\vep(t-r)^{-1}),
\end{split}
\eeq
where in the second to the last inequality we used the fact that $\mu$ is an increasing function.

In conclusion, recall the fact that $p-2=1/p$, we have obtained the lower bound \eqref{lower-2} for $t-r>2$ with a constant $C_0\in(0,1)$ depending only on $f$.
\end{proof}

\subsection{Iteration region}\label{iter-region}
We first identify a region in which the initial estimate
\eqref{lower-2} yields a convenient lower bound for the argument of
$\mu$.  Fix
\[
 \de_0:=\frac{3-p-1/p}4=\frac{3-2\sqrt2}{4}>0.
\]
By \eqref{lowerbd-mu2}, after decreasing $\vep_0$ if necessary (so that $C_{\de_0}\ge\vep_0^{\de_0}$ and $C_0\ge \vep_0^{\de_0}$), we
have
\[
 \mu(\vep\al^{-1})
 \ge C_{\de_0}(\vep\al^{-1})^{\de_0}\ge \vep^{2\de_0}\al^{-\de_0}
\]
whenever $0<\vep\le\vep_0$ and $\al\ge\vep^{-1}$.  Hence, in the
null coordinates $(\al,\be)$, the initial lower bound
\eqref{lower-2} gives
\begin{align}
 u(\al,\be)
 &\ge C_{0}\vep^{p+2\de_0}
 \be^{-1}\al^{-1/p-\de_0} \nonumber\\
 &\ge \vep^{p+3\de_0}\be^{-1}\al^{-1/p-\de_0}\nonumber\\
 &\ge \vep^{p+3\de_0-2}\be^{1-1/p-\de_0}\vep^2\be^{-2}\nonumber\\
 &\ge  \vep^2\be^{-2}\ , \forall \al\ge \vep^{-1} .\label{lower-3}
\end{align}
Here we used
$1-1/p-\de_0=p+3\de_0-2>0
$ and
$\be\ge\al\ge\vep^{-1}$.  Since $\mu$ is increasing, it follows that
\[
 \mu(u(\al,\be))
 \ge \mu( \vep^2\be^{-2})\ .
\]
Substituting this estimate and \eqref{lower-2} into the iteration
formula \eqref{iter-lb}, and restricting the inner integration to
$\vep^{-1}\le\al\le t-r$, we obtain
\beeq\label{u-1}
u(t,r)\ge \frac{C_0^p \vep^{p^2}}{8r}\int_{t-r}^{t+r}\int_{\vep^{-1}}^{t-r}\be^{-p}\al^{-1}\mu^p(\vep\al^{-1})\mu(\vep^2\be^{-2})(\be-\al)d\al d\be\ .
\eneq

To simplify the integral in \eqref{u-1}, we first handle $(\be-\al)$ so that it  behaves like $\be$, by imposing $\al\le\frac{m_{0}}{m_{1}}(t-r)$, for $m_0=\vep^{-1}$ and $m_1>m_0$ to be determined later. In this case, when $\be\ge t-r$, we have $$\be-\al\ge \be-\frac{m_{0}}{m_{1}}(t-r)\ge \frac{m_1-m_0}{m_1}\be\ .$$  
This naturally leads to the following choice of region, that is for $t-r>m_1$, 
\begin{align}
u(t,r)&\ge \frac{C_0^p\vep^{p^2}}{8r}\frac{m_1-m_0}{m_1}\int_{t-r}^{t+r}\int_{m_{0}}^{\frac{m_{0}}{m_1}(t-r)}
\be^{1-p}\al^{-1}\mu^p(\vep\al^{-1})\mu(\vep^2\be^{-2})d\al d\be
\nonumber\\
&\ge \vep^{p^2}\frac{ C_0^p(m_1-m_0)}{8m_1 r}\int_{t-r}^{t+r}\be^{1-p}\mu(\vep^2\be^{-2})d\be
\int_{m_{0}}^{\frac{m_{0}}{m_1}(t-r)}
\al^{-1}\mu^p(\vep\al^{-1})d\al  \label{u-2}
\end{align}
At first, by a change of variables, the integral over $\al$ can be bounded from below in terms of the critical function $A$ introduced in Section \ref{pre}. Indeed, since $m_0=\vep^{-1}$,  we see that
$$
  \int_{m_{0}}^{\frac{m_{0}(t-r)}{m_1}}
\frac{\mu^p(\vep\al^{-1})}{\al}d\al
    =\int_{1}^{\frac{t-r}{m_1}}
\frac{\mu^p(\vep^2 \al^{-1})}{\al}d\al
    = A_{\vep^2}\left(\frac{t-r}{m_1}\right)\ge \frac{1}{2}A_{\vep^2}\left(\left(\frac{t-r}{m_1}\right)^2\right)\,,
$$
where the last inequality follows from Lemma \ref{lem-A} (3).

We now proceed to compute the integral of $\be$. 
Recall that,
 \eqref{lowerbd-mu} with $\de_0=1/2$ implies that for sufficiently small $\vep$, we have
\[
\mu\left(\vep^2\left(2(t-r)\right)^{-2}\right)\ge \frac12\mu(\vep^2(t-r)^{-2})\ .
\]
On the one hand, 
if $t\le 3r$, namely, $t+r\ge2(t-r)$, we derive 
\begin{eqnarray*}
    \frac{1}{2r}\int_{t-r}^{t+r}\be^{1-p}\mu(\vep^2\be^{-2})d\be &\ge& \frac{1}{2r}\int_{t-r}^{2(t-r)} \be^{1-p}\mu(\vep^2\be^{-2})d\be \\
    &\ge& \frac{1}{2r}(t-r)\left(2(t-r)\right)^{1-p}\mu\left(\vep^2\left(2(t-r)\right)^{-2}\right)\\
    &\ge&2^{-p-1}(t+r)^{-1}(t-r)^{2-p}\mu(\vep^2(t-r)^{-2})\, .
\end{eqnarray*}
On the other hand, for $t\ge 3r$, that is $t-r\le t+r\le2(t-r)$, it is similar to obtain that
\begin{eqnarray*}
    \frac{1}{2r}\int_{t-r}^{t+r}\be^{1-p}\mu(\vep^2\be^{-2})d\be &\ge& (t+r)^{1-p}\mu\left(\vep^2(t+r)^{-2}\right)\\
    &\ge&2^{1-p}(t+r)^{-1}(t-r)^{2-p}\mu(\vep^2(t-r)^{-2})\,.
\end{eqnarray*}

As $2-p=-\frac{1}{p}$,  we conclude that  for $t-r>m_1$, 
\begin{equation}\label{lowerbound-1}
    u(t,r)\ge C_0^p\frac{m_1-m_0}{2^{p+4} m_1}\vep^{p^2}(t+r)^{-1}(t-r)^{-\frac1p}\mu(\vep^2(t-r)^{-2})A_{\vep^2}\left(\left(\frac{t-r}{m_1}\right)^2\right)\,.
\end{equation}

\subsection{Iteration procedure}\label{iter-proc}
In this section, we complete the estimate of the upper bound of the lifespan. 
In light of the iteration region and formulation in Section \ref{iter-region}, 
 we establish the following lower bound of the solution by induction.

\begin{prop}\label{H-iterationform}
Let $m_j$ be a sequence of strictly increasing numbers with $m_0=\vep^{-1}$ and $m_j<4/\vep$. 
  Under the assumption of Theorem \ref{blowup-thm}, 
  suppose that for some $j\ge 1$, the solution $u(t,r)$ has the following lower bound estimate:
\beeq\label{iter-lowerbdform}
u(t,r)\ge C_j\vep^{a_j} (t+r)^{-1}(t-r)^{-\frac1p}\mu(\vep^2(t-r)^{-2}) A_{\vep^2}^{b_j}\left(\left(\frac{t-r}{m_{j}}\right)^2\right)\ ,
\eneq
for $t-r> m_{j}$, with some constants
$C_j,a_j,b_j>0$.
Then for any $\de\in(0,1)$, there exists $\vep_0>0$, which is independent of $j\in \N_+$, such that for $0<\vep<\vep_0$, $t-r> m_{j+1}$, \eqref{iter-lowerbdform} holds for $j+1$ with
\beeq\label{iter-lowerbdform-2}
a_{j+1}=a_jp+\de,\ 
b_{j+1}=b_jp+1,\ 
C_{j+1}=
\frac{m_{j+1}-m_j}
{2^{p+4+2\de}(b_jp+1)m_{j+1}}C_j^p.
\eneq
\end{prop}

\begin{proof}

By  \eqref{iter-lb}, \eqref{lower-3} and \eqref{iter-lowerbdform}, for $t-r>m_{j+1}$, we have
\begin{align*}
u(t,r)
&\ge \frac{C_j^p\vep^{a_jp}}{8r}
 \int_{t-r}^{t+r}
 \int_{m_j}^{\frac{m_j}{m_{j+1}}(t-r)}
 \be^{-p}\al^{-1}\mu^p(\vep^2\al^{-2})
 A_{\vep^2}^{b_jp}\left(\left(\frac{\al}{m_j}\right)^2\right)\\
&\qquad\qquad {}
 \times\mu(\vep^2\be^{-2})(\be-\al)\,d\al\,d\be\\
&\ge
 \frac{C_j^p\vep^{a_jp}(m_{j+1}-m_j)}
 {8m_{j+1}r}
 \int_{t-r}^{t+r}
 \be^{1-p}\mu(\vep^2\be^{-2})\,d\be\\
&\qquad\qquad {}
 \times\int_{m_j}^{\frac{m_j}{m_{j+1}}(t-r)}
 \al^{-1}\mu^p(\vep^2\al^{-2})
 A_{\vep^2}^{b_jp}\left(\left(\frac{\al}{m_j}\right)^2\right)
 \,d\al .
\end{align*}

The integral of $\be$ is independent of $j$, and we have already treated it in the proof of \eqref{lowerbound-1}, that is, we have
\beeq \label{integral-be}   \frac{1}{2r}\int_{t-r}^{t+r}\be^{1-p}\mu(\vep^2\be^{-2})d\be \ge 2^{-p-1}(t+r)^{-1}(t-r)^{2-p}\mu(\vep^2(t-r)^{-2})\, .
\eneq

Turning to the integral of $\al$, which plays the key role in the argument. Thanks to \eqref{lowerbd-mu} in Lemma \ref{lem-la}, for any $\de>0$, there exists $\vep_0$ such that for any $\vep<\vep_0$, $\al>m_{j}$ and $j\in\N$,
\beeq\label{loss-de}
\frac{\mu(\vep^2\al^{-2})}{\mu(\vep^2(\al/m_{j})^{-2})}\ge m_{j}^{-\de/p}\ge (\vep/4)^{\de/p}\,.
\eneq
Thus, it implies that
\begin{align*}
    &\int_{m_{j}}^{\frac{m_{j}}{m_{j+1}}(t-r)}\al^{-1}\mu^p(\vep^2\al^{-2})A_{\vep^2}^{b_jp}\left(\left(\frac{\al}{m_{j}}\right)^2\right)d\al\\
    \ge&(\vep/4)^{\de }\int_{m_{j}}^{\frac{m_{j}}{m_{j+1}}(t-r)} \al^{-1}\mu^p(\vep^2(\al/m_{j})^{-2})A_{\vep^2}^{b_jp}\left(\left(\frac{\al}{m_{j}}\right)^2\right)d\al\\
    =&\hf(\vep/4)^{ \de }\int_{1}^{\frac{t-r}{m_{j+1}}}
    A_{\vep^2}^{b_jp}\left(\al^2\right)\pa_\al (A_{\vep^2}(\al^2))d\al =\frac{(\vep/4 )^{\de }}{2(b_jp+1)}A_{\vep^2}^{b_jp+1}\left(\left(\frac{t-r}{m_{j+1}}\right)^2\right)\,.
\end{align*}

Combined with \eqref{integral-be} for the integral of $\be$, we obtain that for any $\de>0$, there exists $\vep_0$ such that for $t-r>m_{j+1}$
\beeq
u(t,r)\ge C_{j+1}\vep^{a_jp+ \de }(t+r)^{-1}(t-r)^{-\frac1p}\mu(\vep^2(t-r)^{-2}) A_{\vep^2}^{b_jp+1}\left(\left(\frac{t-r}{m_{j+1}}\right)^2\right)\ ,
\eneq
where $C_{j+1}=\frac{m_{j+1}-m_{j}}{2^{p+4+2\de }(b_jp+1)m_{j+1}}C_j^{p}$.
\end{proof}

Let us now complete the proof of Theorem \ref{blowup-thm} by estimating the lifespan of the solution. To this end, we choose a sequence of $m_j$ such that
\beeq\label{m-j}
m_j=\frac{4}{\vep}\left(1-2^{-j}\right),\ j\in\N_+.
\eneq

In light of Proposition \ref{H-iterationform} and \eqref{lowerbound-1}, by induction, we deduce that
for $j\in\N_+$,
 and $t-r> m_{j+1}$, we have \eqref{iter-lowerbdform}, with
\beq
\begin{split}
a_{j+1}&=a_jp+\de=p^{j+2}+\frac{p^j-1}{p-1}\de\,,\\
b_{j+1}&=b_jp+1=p^{j}+p^{j-1}+\cdots+1=\frac{p^{j+1}-1}{p-1}\le p^{j+1}\ ,\\
C_1&=C_0^p\frac{m_1-m_0}{2^{p+4} m_1}=\frac{1}{2^{p+5}}C_0^{p}\ ,\\
C_{j+1}&=\frac{ m_{ j+1 }-m_{ j}}{
2^{p+4+2\de }
(b_jp+1)m_{ j+1 }}C_j^{p}\ge\frac{1}{2^{9}(2p)^{j+1} }C_j^{p}\, .
\end{split}
\eeq
Thus, the estimate of $C_j$ yields
\beq
\begin{split}
\log C_{j+1}&\ge -9 \log 2-(j+1)\log 2p+ p\log C_{j}\\
&\ge p^{j}\left(\log C_1-(\log 2p)\sum_{l=1}^{\infty}\frac{l+1}{p^l}-(9\log 2)\sum_{l=1}^{\infty}\frac{1}{p^l}\right)=:\ti{C}p^j\ .
\end{split}
\eeq 
Finally, for $t-r>4\vep^{-1}>m_{j+1}$, we have 
\beq
\begin{split}
u(t,r)
&\ge
C_{j+1}\vep^{a_{j+1} }(t+r)^{-1}(t-r)^{-\frac1p}\mu(\vep^2(t-r)^{-2}) A_{\vep^2}^{b_{j+1}}\left(\left(\frac{t-r}{m_{j+1}}\right)^2\right)\\
&\ge \exp\left[{p^j}\left(\ti{C}+\left(p^2+\de\frac{1-p^{-j}}{p-1}\right)\log\vep+\frac{p-p^{-j}}{p-1}\log A_{\vep^2}\left(\left(\frac{\vep(t-r)}{4}\right)^2\right)\right)\right]\\
&\ \ \times (t+r)^{-1}(t-r)^{-\frac1p}\mu(\vep^2(t-r)^{-2})\ ,
\end{split}
\eeq
for any $j\ge 1$. 
According to $\frac{p-\frac1{p^j}}{p-1}\to \frac{p}{p-1}$ 
and $p^2+\de\frac{1-p^{-j}}{p-1}\to p^2+\frac{\de}{p-1}$, 
we see that $u(t,r)$ explodes as $j\to\infty$, provided that 
$$\ti{C}+\left(p^2+\frac{\de}{p-1}\right)\log\vep+\frac{p}{p-1}\log A_{\vep^2}\left(\left(\frac{\vep(t-r)}{4}\right)^2\right)>0.$$
In other words, for such $(t,r)$, the solution $u(t,r)$ cannot exist, which implies that the solution blows up in finite time.

In particular, for sufficiently small $\vep$, with $r=0$, the infimum of such $t$ gives an upper bound for the lifespan $T(\vep)$. That is, we  have the following implicit upper bound for the lifespan $T(\vep)$:
\beeq\label{lifespan-up}
\ti{C}+\frac{p}{p-1}\log A_{\vep^2}\left(\left(\frac{\vep T(\vep)}{4}\right)^2\right)\le \left(p^2+\frac{\de}{p-1}\right)\log\frac1{\vep}\,.
\eneq

The change of variables $\vep\la=\ti{\la}$ and Lemma \ref{lem-A} yield
$$\hf A_{\vep^2}((\vep T(\vep))^2)\le A_{\vep^2}(\vep T(\vep))=A_{\vep}(T(\vep))-A_{\vep}(1/\vep)\,,$$
which shows by the upper bound \eqref{lifespan-up} that
$$A_{\vep}(T(\vep))-A_{\vep}(1/\vep)\les \vep^{-p(p-1)-\de/p}\,.$$
Recall that by the expression of $A_{\vep}(t)$, 
$$A_{\vep}(1/\vep)\les\int_{\vep}^1\frac1{\la}d\la=\log 1/\vep\,,$$
we obtain the desired upper bound of the lifespan:
\beeq\label{lifespan-up2}
A_{\vep}(T(\vep))\les \vep^{-p(p-1)-\de/p}\,.\eneq

\begin{remark}\label{mu=1}
    When $\mu\equiv1$, that is, for the classical John-Strauss nonlinearity $|u|^{1+\sqrt{2}}$, we could choose $\de=0$ for the lower bound in \eqref{loss-de}. In this case, $A_{\vep^2}(t)= \log t$  and \eqref{lifespan-up} reduces to 
    $$
    \ti{C}+\frac{p}{p-1}\log \left(2\log\frac{\vep T(\vep)}4\right)\le  p^2 \log\frac1{\vep}\, ,
    $$
    which implies that, for some $C>0$ independent of $\vep$,
    $$\log\frac{\vep T(\vep)}4\le C{\vep}^{-p(p-1)}\, .$$
    For sufficiently small $\vep$, we get the upper bound of the lifespan
    $$T(\vep)\le 4\vep^{-1}\exp(C{\vep}^{-p(p-1)})
    \le \exp(2C{\vep}^{-p(p-1)})   
    \,,$$
    which is  consistent with the well known critical upper bound.
\end{remark}

\section{The lower bound estimate of the lifespan}\label{lower-esti-life}
In this section, to compare with the upper bound \eqref{lifespan-up2}, we establish the lower bound of the lifespan $T(\vep)$, \eqref{eq-life}, for the solution to \eqref{blow-eq} with small initial data. We remark here that the result applies also for general data, with $u(0,x)=\vep \phi$ instead of $u(0,x)=0$.

We establish the lower bound by using an argument similar to the proof of global well-posedness in Section \ref{proof-global}.
Suppose that the initial data have compact support in $B_1(O)$, then by Lemma \ref{linear-esti} there exists a constant $B_0$ depending only on $f$ such that
\beeq\label{inside-sone}
|u^0(t,r)|\le B_0\vep\<t+r\>^{-1}\chi_{|t-r|\le1}:=\frac 14 \Psi_{in}\,.
\eneq
Similar to
the proof of Lemma \ref{bound-lem}, we derive the 
estimate 
as follows
\beeq\label{psi-in}
\begin{split}
   | LF(u^0)(t,r)|
   &\le | LF(\Psi_{in})(t,r)|
   \\
   &\le \frac{(4 B_0\vep)^p}{8r}\int_{|t-r|}^{t+r}\int^{t-r}_{-\be} 
\<\be\>^{-p}\mu(4B_0\vep\<\be\>^{-1})\chi_{   |\al|\le1}
(\be-\al)d\al d\be\\
    &\le \frac{(4 B_0\vep)^p}{8r}\mu(4B_0\vep\<t-r\>^{-1})\int_{|t-r|}^{t+r}\<\be\>^{1-p}d\be\\
    &\le B_1\vep^{p}\<t+r\>^{-1}\<t-r\>^{-\frac1p}\mu(\vep\<t-r\>^{-1})\,,
\end{split}
\eneq
for some constant $B_1>B_0$ independent of $\vep\le 1$. 

Thus, we establish the ansatz of the upper bound of the solution 
\beeq
\begin{split}
    |u(t,r)|&\le 4B_0\vep\<t+r\>^{-1}\chi_{|t-r|\le1}+4B_1\vep^{p}\<t+r\>^{-1}\<t-r\>^{-\frac1p}\mu(\vep\<t-r\>^{-1})\chi_{t-r>1}\\
    &:=\Psi_{in}+\Psi_{out}\ .
\end{split}
\eneq
 As in Lemma \ref{bound-lem}, the local solution $u(t,r)$ could be extended to the interval $t\in [0,T]$, provided that the bootstrap assumption holds on $[0,T]$. In particular, it reduces to the following estimate:
\beeq\label{eq-bootstrap}
LF(\Psi_{in}+\Psi_{out})(t,r)\le \hf(\Psi_{in}+\Psi_{out})(t,r)\,, \ \forall (t,r)\in[0,T]\times\R_+\,.
\eneq

At first, by \eqref{psi-in}, 
it is clear that there exists $\vep_0<1$ such that for any $\vep<\vep_0$, 
\beeq\label{vep-1}
B_1\vep^{p-1}
\<t-r\>^{-\frac1p}\mu(\vep\<t-r\>^{-1})
\le
B_1\vep_0^{p-1}
\mu(\vep_0)
\le   B_0\wedge \frac{1}{4 } \,,
\eneq
and so
$\Psi_{out}\le \vep\<t+r\>^{-1}$, 
$$LF(\Psi_{in})(t,r)\le \frac14(\Psi_{in}+\Psi_{out})(t,r)\,.$$

On the other hand, as $\Psi_{out}$ is supported inside the light cone, the contribution to the solution has the same property in the support. In conclusion, it follows that, for $t-r>1$, we have
\beeq
\begin{split}
    LF(\Psi_{out})(t,r)&\le \frac{(4B_1\vep^p)^p}{8r}\int_{ t-r }^{t+r}\int_{1}^{t-r}\<\be\>^{-p}\<\al\>^{-1}\mu^p(\vep\<\al\>^{-1})\mu(\vep\<\be\>^{-1})(\be-\al)d\al d\be\\
    &\le \frac{(4B_1\vep^p)^p}{8r}\mu(\vep\<t-r\>^{-1})\int_{ t-r }^{t+r}\<\be\>^{1-p}d\be \int_{1}^{\be}\<\al\>^{-1}\mu^p(\vep\<\al\>^{-1})d\al\\
    &\le (4B_1\vep^p)^p\mu(\vep\<t-r\>^{-1})\frac{1}{8r}\int_{ t-r }^{t+r}\<\be\>^{1-p}A_{\vep}(\<\be\>)d\be\\
    &\le  (4B_1\vep^p)^p \mu(\vep\<t-r\>^{-1})A_{\vep}(\<t+r\>)
    \frac{1}{8r}\int_{ t-r }^{t+r}\<\be\>^{1-p} d\be
    \\
    &\le B_2\vep^{p^2}\<t+r\>^{-1}\<t-r\>^{-\frac1p}\mu(\vep\<t-r\>^{-1})A_{\vep}(1+2t)
\end{split}
\eneq
for some constant $B_2$ independent of $\vep$, where the last inequality follows from
 the last two 
 inequalities in the 
 proof of Lemma \ref{bound-lem}, as well as the trivial inequality $\<t+r\>=2+t+r\le 1+2t$. 

Therefore, we can close the bootstrap argument of \eqref{eq-bootstrap} on any time interval $[0,T]$, provided that
$$B_2\vep^{p^2} A_{\vep}(1+2t)\le B_1\vep^p\ .$$
Consequently, by the definition of the lifespan, $T(\vep)$ must satisfy
$$B_2\vep^{p^2} A_{\vep}(1+2T(\vep))\ge B_1\vep^p\,.$$

Recall that $T(\vep)\ge 1$ for small $\vep$,   
 the proof of Lemma \ref{lem-A} (2) with $\theta=1/3$ gives us
 $$A_{\vep}(1+2T(\vep))\le A_{\vep}(3T(\vep))\le 3 A_{\vep}( T(\vep))\ ,$$
and thus the desired lower bound of the lifespan in
\eqref{eq-life},
with $c=B_1/(3B_2)$.


\bibliography{bib.bib}
 \end{document}